\documentclass[11pt]{article}
\usepackage{mathrsfs}
\usepackage{amsmath,amssymb}
\def\P{\mathcal{P}}
\def\PP{\mathscr{P}}
\def\A{\mathcal{A}}
\def\Y{\mathscr{G}}
\def\E{\mathscr{E}}
\def\R{\mathscr{R}}
\def\C{\mathbb{C}}
\def\D{\mathscr{D}}
\def\T{\mathbb{T}}
\newtheorem{theorem}{\hspace*{\parindent}Theorem}
\newtheorem{lemma}{\hspace*{\parindent}Lemma}
\newtheorem{corollary}{\hspace*{\parindent}Corollary}
\title{Four-point distortion theorem for complex polynomials}
\author{V.N.\,Dubinin\\[15pt]\small{Institute of
Applied Mathematics, 7 Radio Street, Vladivostok, 690041, Russia}\\[5pt]\small{E-mail: dubinin@iam.dvo.ru}
\\[15pt]
Dedicated to the memory of Professor Promarz M. Tamrazov}

\date{}
\begin{document}
\maketitle

\begin{abstract}
We prove a theorem on distortion of cross ratio of four points under the mapping
effected by a complex polynomial with restricted  critical
values.  Its corollaries include inequalities involving the absolute
value and certain coefficients of a polynomial.  In particular, an
exact lower bound is established for  maximal moduli of critical values of polynomials $P$
of degree $n$ normalized by $P(0)=0$, $P'(0)\ne0$.
\end{abstract}

\bigskip


\bigskip


\bigskip

\section{Introduction}
It is a common knowledge that certain inequalities for complex
polynomials may be derived using the theory of univalent functions
\cite[Chapter~2]{DubininUMN}.  Solutions of many extremal problems
in that theory hinge on metric and conformal properties of
multiply connected domains (see, for instance,
\cite{Goluzin}-\cite{Tamrazov2}).  We have  recently amplified
this approach by extending the well-known extremal problems of
Gr\"{o}tzsch and Teichm\"{u}ller for moduli of doubly connected
planar domains to domains located on certain Riemann surfaces
\cite{DubininMZ}.  In this note we consider some applications of
the solution of such extended Teichm\"{u}ller problem to
inequalities for polynomials with restricted critical values.  Let
us remind that a critical value of the polynomial $P$ is its value
$P(\zeta)$ at a point $\zeta$, where $P'(\zeta)=0$.  Denote by $\P_n$, $n\geq{2}$, the class of all complex polynomials of degree $n$ whose critical values do not exceed unity in absolute value.  This class is quite rare in the literature (cf. \cite{MMR}-\cite{RS}).   Note that under more traditional restriction $P(z)=0\Rightarrow|z|<1$ all critical points of $P$ also lie in the unit disk which implies
$$
\max\limits_{P'(\zeta)=0}|P(\zeta)|\leq\max\limits_{|z|<1}|P(z)|.
$$
In this case an upper bound imposed on the uniform norm of a polynomial
leads to its membership in $\P_n$ .   In section~2 we will prove a distortion theorem for
cross ratio of four points under a mapping involving an arbitrary polynomial $P$ from
$\P_n$ and the Chebyshev polynomial of the first kind $w=T_n(z)=2^{n-1}z^n+\cdots$  We
will need a more comprehensive description of the Chebyshev polynomial $T_n(z)$. It can be defined in terms of conformal maps as the composition of the inverse Zhukowski map, the power function and the direct Zhukowski map:
$$
T_n(z)=\frac{1}{2}\left(\left(z+\sqrt{z^2-1}\right)^n+\left(z-\sqrt{z^2-1}\right)^n\right),~~~z\in\overline{\C}.
$$
The hyperbolas with foci at $z=\pm1$ passing through the critical points $z=\cos(\pi{k}/n)$, $k=1,\ldots,n-1$, of the polynomial  $T_n(z)$ partition the $z$-plane into $n$ pairwise disjoint domains.  Let $B_1,\ldots,B_n$ denote these domains numbered from right to left. The polynomial $T_n$ effects univalent conformal mapping of $B_1$ onto the domain $D_1$ which is entire $w$-plane cut along the ray $L^-:=[-\infty,-1]$.  The domains $B_2,\ldots,B_{n-1}$ are mapped by this polynomial onto the domains $D_2,\ldots,D_{n-1}$ all of which are copies of the $w$-plane cut along the rays $L^-$ and $L^+:=[1,\infty]$. Finally, $B_n$ is mapped onto $D_n=\overline{\C}\setminus{L^-}$ if $n$ is even or $D_n=\overline{\C}\setminus{L^+}$  if $n$ is odd.
We can construct the Riemann surface  $\R(T_n)$ of the function inverse  to $T_n$ by gluing together the domains $D_k$, $k=1,\ldots,n$, as follows: $D_1$ is glued crosswise to $D_2$ along the sides of the cuts made along the ray $L^-$.  Domain $D_2$ is glued to $D_3$ along the sides of the cuts made along the ray $L^+$, and so on.  Domain $D_{n-1}$ is glued to $D_n$ through $L^-$ if $n$ is even and through $L^+$ is $n$ is odd.  The domain $D_k$ viewed as a subset of the Riemann surface $\R(T_n)$  will be denoted by $\D_k$.  The proof of the main result of this paper (Theorem~\ref{th:main}) hinges on an extremal property of an analogue of the Teichm\"{u}ller ring lying on the surface  $\R(T_n)$. Perhaps, for researchers applying polynomial inequalities more interesting are corollaries of Theorem~1 collected in Section~3.  In particular, Corollary~4 containing a lower bound for the maximal moduli of critical values of a polynomial may be of interest.

\section{The main result}
Let
$$
(a_1,a_2,a_3,a_4)=\frac{a_3-a_1}{a_3-a_2}:\frac{a_4-a_1}{a_4-a_2}
$$
be the anharmonic ratio (or cross ratio) of four ordered distinct  points $a_k$, $k=1,\ldots,4$, lying in the extended complex plane $\overline{\C}$.  For any given polynomial $P$ and a point $z\in\overline{\C}$ denote by $x_P(z)$ the root of the equation $T_n(x)=|P(z)|$ lying on the ray $[\cos(\pi/(2n)),+\infty]$.  The number $\cos(\pi/(2n))$ here is the largest zero of the Chebyshev polynomial $T_n$.

\begin{theorem}\label{th:main}
Suppose $P\in\P_n$. Then for any four distinct points $z_k$, $k=1,\ldots,4$, located on an oriented straight line in ascending order the following inequality holds\emph{:}
\begin{equation}\label{eq:main}
-(z_3,z_1,z_2,z_4)\leq|(x_P(z_3),-x_P(z_1),-x_P(z_2),x_P(z_4))|.
\end{equation}
Equality in \emph{(\ref{eq:main})} is attained, for instance, for  $P=T_n$ and any points $z_k$, $k=1,\ldots,4$, satisfying $-\infty<z_1<z_2<-\cos(\pi/(2n))$, $\cos(\pi/(2n))<z_3<z_4<+\infty$.
\end{theorem}
In order to give a proof of the theorem we will need some definitions and results from \cite{DubininMZ}.
Here and in what follows the Riemann surface  $\R$ is understood as bordered compact  Riemann surface.  We view  it as lying over the sphere $\overline{\C}_w$ and made up of finite number of planar domains with natural definitions of projection, local parameter and neighborhood for points on such surface \cite{CourantHurwitz}. Following \cite{DubininMZ}, denote by $\A_n$, $n\ge{2}$,  the collection of doubly connected domains $\Y$ satisfying the following conditions:

i) the ring $\Y$ is located on  a Riemann surface $\R$ covering every point of the sphere $\overline{\C}_w$  not more than $n$ times;

ii) the complement $\R\!\setminus\!\Y$ consists of two connected components $\E_0$ and $\E_1$ one of which ($\E_0$) contains the entire boundary of $\R$;

iii) any closed Jordan curve in $\R\!\setminus\!\E_0$ that lies over the circle $|w|=\rho$, $1\leq\rho<\infty$, forms an $n$-fold covering of that circle.

By the Teichm\"{u}ller ring on the surface $\R(T_n)$, denoted by
$$
\T_n(s,t,\sigma,\tau),~~n\ge2,~~0\leq{s}<t<\infty, ~~0\leq{\sigma}<\tau\leq\infty,
$$
we will mean the doubly connected domain obtained from $\R(T_n)$ by deleting the interval with projection $[(-1)^ns,(-1)^nt]$ from the sheet $\D_n$ and the interval  with projection  $[\sigma,\tau]$ from the sheet $\D_1$.

  The next claim is true.

\begin{lemma}\label{lm:Teichm}\emph{(\cite[theorem~2]{DubininMZ})}.  The Teichm\"{u}ller ring $\T_n(s,t,\sigma,\tau)$ has maximal modulus among all doubly connected domains $\Y\in\A_n$ such that the projection of one component  $\E_0$ of the complement $\R\!\setminus\!\Y$ connects the circles $|w|=s$ and $|w|=t$ while the projection of the other component  $\E_1$ connects the circles
$|w|=\sigma$ and $|w|=\tau$ \emph{(}$0\leq{s}<t<\infty$, $0\leq\sigma<\tau\leq\infty$\emph{)}.
\end{lemma}
Recall that the modulus of doubly connected domain with respect to the family of curves separating its boundary components is defined to be $(2\pi)^{-1}\log(R_2/R_1)$, where $R_1<R_2$ are the radii of the inner and outer boundary circles of any annulus conformally equivalent  to the given domain \cite{Jenkins}.

\textbf{Proof of Theorem~\ref{th:main}:} The equality case in (\ref{eq:main}) can be verified directly from the definition of the Chebyshev polynomial.  To prove inequality (\ref{eq:main}) it suffices to consider the case when all critical values of the polynomial $P$ lie in the open unit disk while the points $z_k$, $k=1,\ldots,4$, are real and ordered as $z_1<z_2<z_3<z_4$ and $|P(z_1)|\ne|P(z_2)|$, $|P(z_3)|\ne|P(z_4)|$.  Denote by $\PP^{-1}$ the inverse function to the polynomial $P$. It is one-valued and analytic on the Riemann surface $\R(P)$.  Let $\PP:\overline{\C}_z\to\R(P)$ be its inverse. Let doubly connected domain $G$ on the Riemann sphere be defined by the two components $E_0=[z_1,z_2]$ and $E_1=[z_3,z_4]$ of its complement.  Suppose $\Y=\PP(G)$, $\E_0=\PP(E_0)$, $\E_1=\PP(E_1)$. Conformal invariance of the modulus implies $\mod{G}=\mod{\Y}$. We claim that the doubly connected domain $\Y\subset\R:=\R(P)\!\setminus\!\E_0$ belongs to the collection $\A_n$.  Indeed, conditions i) and ii) are obviously satisfied. To verify condition iii)  first note that that the part of $\R(P)$ lying over the domain $|w|>\rho$, $\rho>1$, is itself a Riemann surface forming $n$-fold covering of this domain.  Since $P\in\P_n$ the only branch point of this surface is at infinity and has order $n-1$. If condition iii) is violated then there is more than one boundary curve of that surface lying over the circle $|w|=\rho$ which contradicts the Hurwitz formula.  Hence, $\Y\in\A_n$.  In addition, the projection of $\E_0$ coincides with $P(E_0)$ and is a connected set containing the points $P(z_1)$ and $P(z_2)$.  Similarly,
the projection of $\E_1$ coincides with $P(E_1)$ and is a connected set containing the points $P(z_3)$ and $P(z_4)$. Assume  that $|P(z_1)|<|P(z_2)|$, $|P(z_3)|<|P(z_4)|$. According to Lemma~\ref{lm:Teichm}
$$
\mod{\Y}\leq\!\!\!\!\!\mod{\T_n\big(|P(z_1)|,|P(z_2)|,|P(z_3)|,|P(z_4)|\big)}.
$$
Conformal invariance of the modulus now yields
$$
\mod{\T_n\big(|P(z_1)|,|P(z_2)|,|P(z_3)|,|P(z_4)|\big)}=\!\!\!\!\!\mod\!\tilde{G},
$$
where $\tilde{G}=\overline{\C}_z\!\setminus\!\big([-x_P(z_2),-x_P(z_1)]\cup[x_P(z_3),x_P(z_4)]\big)$. Finally we arrive at
\begin{equation}\label{eq:GtildeG}
\mod\!G\leq\!\!\!\!\!\mod\!\tilde{G}.
\end{equation}
If $|P(z_1)|>|P(z_2)|$ or  $|P(z_3)|>|P(z_4)|$ we still get (\ref{eq:GtildeG}) in a similar manner.
Let $\varphi$ be the linear fractional automorphism of the sphere $\overline{\C}_z$ sending the points $z_1$, $z_2$, $z_4$ to the points $-1$, $0$, $\infty$, respectively, and let $\psi$  be the linear fractional automorphism of the sphere $\overline{\C}_z$ sending the points $-x_P(z_1)$, $-x_P(z_2)$, $x_P(z_4)$ to the points $-1$, $0$, $\infty$, respectively.
In view of (\ref{eq:GtildeG}),
$$
\mod\!\varphi(G)\leq\!\!\!\!\!\mod\!\psi(\tilde{G})
$$
so that
$$
\varphi(z_3)\leq|\psi(x_P(z_3))|.
$$
It is left to notice that
$$
\varphi(z_3)=-(\varphi(z_3),-1,0,\infty)=-(z_3,z_1,z_2,z_4)
$$
and
$$
\psi(x_P(z_3))=|(\varphi(x_P(z_3)),-1,0,\infty)|=|(x_P(z_3),-x_P(z_1),-x_P(z_2),x_P(z_4)|.\hfil \square
$$

\textbf{Remark~1}. The proof of Theorem~\ref{th:main} shows that given any four points $z_k$, $k=1,\ldots,4$, on $\overline{\C}_z$ and a doubly connected domain $G$ separating the pair $z_1$, $z_2$ from the pair $z_3$, $z_4$ we will have
$$
\mod\!G\leq\!\!\!\!\!\mod\!\{\overline{\C}_z\!\setminus\!\big([-x_P(z_1),-x_P(z_2)]\cup[x_P(z_3),x_P(z_4)]\big)\}.
$$

\textbf{Remark~2}. It appears to be an interesting problem to establish an analogue of inequality (\ref{eq:main}) for arbitrary $n$-valent functions in the disk $|z|<1$ such that their critical values do not exceed unity in absolute value (cf. \cite[inequality (1)]{DubKost}).

\section{Corollaries}

In this section we will apply Theorem~\ref{th:main} to get some estimates for polynomials
$$
P(z)=c_0+c_1z+\cdots+c_nz^n
$$
of the class $\P_n$ involving the values of $|P(z)|$ at certain points $z$ or containing the coefficients $c_1,\ldots,c_n$.
\begin{corollary}\label{cr:1}
Suppose $P\in\P_n$, $z_1$ and $z_2$ are arbitrary point of $\C_z$ such that $|P(z_k)|\leq1$, $k=1,2$.  Then
$$
|P(z)|\leq T_n\left(\Big|\frac{2z-z_1-z_2}{z_1-z_2}\Big|\right)
$$
for all points $z$ lying of the rays $z=\frac{1}{2}((z_2-z_1)t+z_2+z_1)$, $t\leq-1$ and $t\geq1$.
\end{corollary}
\textbf{Proof.}
Suppose $P\in\P_n$ satisfies $|P(\pm1)|\leq1$.  Applying theorem~\ref{th:main} to this polynomial and points $\infty$, $-1$, $1$, $t$ ($t>1$) we arrive at
$$
\frac{2}{t-1}\leq\frac{2}{|x_P(t)-1|}.
$$
This implies $t\geq{x_P(t)}$ so that
$$
|P(t)|\leq\T_n(t)~\mathrm{for}~t>1.
$$
Considering similarly the polynomial $P(-t)$ and the same points we get
$$
|P(-t)|\leq\T_n(t),
$$
leading to
\begin{equation}\label{eq:PT}
|P(t)|\leq\T_n(|t|)~\mathrm{for all}~|t|>1.
\end{equation}
General case reduces to the one just proved by the change of variable  $z=\frac{1}{2}((z_2-z_1)t+z_2+z_1)$.
$\square$
Note that inequality (\ref{eq:PT}) has the same form as the classical Chebyshev inequality (see, for instance, \cite[p.235]{BorweinErdelyi}).
\begin{corollary}\label{cr:2}
Suppose $P(z)=c_0+c_1z+\cdots+c_nz^n$, $c_n\ne0$, is an arbitrary  polynomial with complex coefficients and $z_1$, $z_2$ are any points in $\C_z$  satisfying
\begin{equation}\label{eq:points}
|z_1-z_2|>\big(2^{2n-1}/|c_n|\big)^{1/n}.
\end{equation}
The either $|P(z_1)|>1$ or $|P(z_2)|>1$ or there exists a critical value of $P$ with absolute value strictly greater than one.
\end{corollary}
\textbf{Proof.}
 Suppose the conclusion of Corollary~\ref{cr:2} is wrong. Then all hypotheses of Corollary~\ref{cr:1} are satisfied.  However, the conclusion of Corollary~\ref{cr:1} contradicts (\ref{eq:points}) if $|z|$ is sufficiently large.
$\square$

\begin{corollary}\label{cr:3}
If the polynomial $P(z)=c_1z+\cdots+c_nz^n$, $c_1\ne0$, belongs to $\P_n$ then
\begin{equation}\label{eq:PbiggerT}
|P(z)|\geq T_n\Big(\frac{1}{n}(\sin\frac{\pi}{2n})|c_1z|-\cos\frac{\pi}{2n}\Big)
\end{equation}
for all points whose absolute values are greater than or equal to $2n\cot(\pi/(2n))/|c_1|$.  Equality is attained for
$P(z)=T_n\big(z-\cos(\pi/(2n))\big)$ at the points $z\geq2\cos(\pi/(2n))$.
\end{corollary}
\textbf{Proof.}
Fix $r>0$ and $z$ satisfying the condition of the corollary.  Set $z_1=-rz/|z|$, $z_2=0$, $z_3=z$ and $z_4=\infty$. Then $x_P(z_1)=T_n^{-1}\big(|P(-rz/|z|)|\big)$ represents the value of the inverse function $T_n^{-1}$ lying on the ray $[\cos(\pi/(2n)),+\infty]$, $x_P(z_2)=\cos(\pi/(2n))$, $x_P(z_3):=x_P(z)$ and $x_P(z_4)=\infty$.  Inequality (\ref{eq:main}) then takes the form
$$
\frac{|z|}{r}\leq\left|\frac{-\cos(\pi/(2n))-x_P(z)}{-\cos(\pi/(2n))+T_n^{-1}\big(|P(-rz/|z|)|\big)}\right|.
$$
Multiplying both sides by $r$ and taking limit as $r\to0$ we obtain
$$
|c_1z|\leq\big|T_n'(\cos(\pi/(2n)))\big|\big|\cos(\pi/(2n))+x_P(z)\big|\leq\frac{n}{\sin\frac{\pi}{2n}}\left(\cos\frac{\pi}{2n}+T_n^{-1}\big(|P(z)|\big)\right),
$$
where the value of the inverse function is again chosen to belong to the ray $[\cos(\pi/(2n)),+\infty]$.  Since $T_n$ is increasing on that ray we have the required inequality.  Equality case can be verified directly.
$\square$
\begin{corollary}\label{cr:4}
For any polynomial $P(z)=c_1z+\cdots+c_nz^n$, $c_1\ne0$, $c_n\ne0$, there exits a critical value $P(\zeta)$ \emph{(}$P'(\zeta)=0$\emph{)} such that
\begin{equation}\label{eq:critical}
|P(\zeta)|\ge2\left(\frac{1}{n}\sin\frac{\pi}{2n}\right)^{\frac{n}{n-1}}\left|\frac{c_1^n}{c_n}\right|^{\frac{1}{n-1}}.
\end{equation}
The constant on the right hand side of \emph{(\ref{eq:critical})} cannot be made bigger.
\end{corollary}
\textbf{Proof.} Put $M=\max\{|P(\zeta)|:~P'(\zeta)=0\}$.  The the polynomial $P/M$ belongs to $\P_n$. Substituting this polynomial into  (\ref{eq:PbiggerT}) and comparing coefficients at $|z|^n$ in (\ref{eq:PbiggerT}) we are led to the estimate (\ref{eq:critical}) with $|P(\zeta)|=M$.  Equality in (\ref{eq:critical}) is attained for the polynomial $T_n(z-\cos(\pi/(2n)))$ whose critical values are unimodal.
$\square$
It seems interesting to compare inequality (\ref{eq:critical}) with the upper bound for the critical values which we obtained earlier in \cite{DubSbornik}:
$$
|P(\zeta)|\leq(n-1)\left(\frac{1}{n}\right)^\frac{n}{n-1}\left|\frac{c_1^n}{c_n}\right|^\frac{1}{n-1}.
$$
Here $P$ represents a polynomial from corollary~\ref{cr:4} and $\zeta$ is one of its critical points.  Equality is attained if $P(z)=c_1z+c_nz^n$.
\begin{corollary}\label{cr:5}
If a polynomial $P(z)=\sum_{k=0}^{n}c_kz^k$ with real coefficients $c_k$, $k=0,1,\ldots,n$,  belongs to $\P_n$ then
\begin{equation}\label{eq:realcoef}
c_nc_{n-2}+n2^{-\frac{2}{n}}c_n^{2-\frac{2}{n}}\ge\frac{n-1}{2n}c_{n-1}^2.
\end{equation}
Equality occurs for the Chebyshev polynomial $T_n$.
\end{corollary}
\textbf{Proof.}
There is no loss of generality in assuming $c_n>0$.  The function $T_n^{-1}(P(z))$ is analytic in the neighborhood of infinity.   Let $f(z)$ denote its branch which is positive for positive $z$.  We will show first that for sufficiently large positive $x$ the following holds
\begin{equation}\label{eq:positivx}
\left|\frac{4x^2f'(x)f'(-x)}{(f(x)-f(-x))^2}\right|\leq1.
\end{equation}
To this end apply Theorem~\ref{th:main} to polynomial $P$ and points
$z_1=-x-\Delta{x}$, $z_2=-x$, $z_3=x$ and $z_4=x+\Delta{x}$, $\Delta{x}>0$. For these points:
\begin{equation*}
\begin{split}
&x_P(z_1)=T_n^{-1}\big(|P(-x-\Delta{x})|\big)=-f(-x-\Delta{x}),
\\
&x_P(z_2)=T_n^{-1}\big(|P(-x)|\big)=-f(-x),
\\
&x_P(z_3)=T_n^{-1}\big(|P(x)|\big)=f(x),
\\
&x_P(z_4)=T_n^{-1}\big(|P(x+\Delta{x})|\big)=f(x+\Delta{x}),
\end{split}
\end{equation*}
where $T_n^{-1}(\cdot)$ is the value of the inverse function  $T_n^{-1}$ lying on the ray $[\cos(\pi/(2n)),+\infty]$.  Inequality (\ref{eq:main}) yields
$$
\frac{4x(x+\Delta{x})}{(\Delta{x})^2}\leq\left|\frac{(f(-x)-f(x))(f(x+\Delta{x})-f(-x-\Delta{x}))}{(f(-x)-f(-x-\Delta{x})(f(x+\Delta{x})-f(x)))}\right|.
$$
Multiplying both sides by $(\Delta{x})^2$ and passing to the limit as $\Delta{x}\to0$ we arrive at inequality (\ref{eq:positivx}).

Further, since
$$
y=T_n(x)=2^{n-1}x^n-n2^{n-3}x^{n-2}+o(x^{n-2}),~~x\to+\infty,
$$
 we will have for positive values of the root
$$
\sqrt[n]{y}=2^{\frac{n-1}{n}}x(1-\frac{1}{4}x^{-2}+o(x^{-2})),~~x\to+\infty.
$$
Hence,
$$
T_n^{-1}(y)=x=2^{\frac{1-n}{n}}\sqrt[n]{y}+2^{\frac{-n-1}{n}}\frac{1}{\sqrt[n]{y}} +o\left(\frac{1}{\sqrt[n]{y}}\right),~~y\to+\infty.
$$
On the other hand,
\begin{multline*}
\sqrt[n]{P(x)}=\sqrt[n]{c_n}x\left(1+\frac{c_{n-1}}{c_n}\frac{1}{x}+\frac{c_{n-2}}{c_n}\frac{1}{x^2}+o\left(\frac{1}{x^2}\right)\right)^{\frac{1}{n}}
\\
=\sqrt[n]{c_n}\left\{x+\frac{1}{n}\frac{c_{n-1}}{c_n}+\left[\frac{1}{n}\frac{c_{n-2}}{c_n}-\frac{n-1}{2n^2}\frac{c_{n-1}^2}{c_n^2}\right]\frac{1}{x}+o\left(\frac{1}{x}\right)\right\}, ~~x\to+\infty,
\end{multline*}
so that for positive $x$
$$
f(x)=2^{\frac{1-n}{n}}\sqrt[n]{c_n}\left(x+\frac{1}{n}\frac{c_{n-1}}{c_n}+d\frac{1}{x}+o\big(\frac{1}{x}\big)\right),~~x\to+\infty,
$$
where
$$
d=\frac{1}{n}\frac{c_{n-2}}{c_n}-\frac{n-1}{2n^2}\frac{c_{n-1}^2}{c_n^2}+2^{-\frac{2}{n}}\frac{1}{\sqrt[n]{c_n^2}}.
$$
Due to analyticity this expansion holds in some neighborhood of infinity and, in particular, for some  negative $x$.  Substituting the above expansion of $f$ into (\ref{eq:positivx}) leads to inequality $d\ge0$ which is equivalent to (\ref{eq:realcoef}). Equality can be verified by substituting the coefficients of the Chebyshev polynomial into (\ref{eq:realcoef}).
$\square$

\section{Acknowledgements}
This work has been supported by the Russian Foundation for Basic Research (grant no.11-01-0038) and Far Eastern Branch of the Russian Academy of Sciences (grant no. 12-I-OMH-02).


\begin{thebibliography}{99}

\bibitem{DubininUMN} V.N.Dubinin, \emph{Methods of geometric function theory in classical and modern problems for
polynomials}, Russian Math. Surveys, 67(4)(2012), pp.599--684.

\bibitem{Goluzin} G.M.\,Goluzin, \emph{Geometric theory of functions of
a complex variable}, Translations of Mathematics Monographs,
volume 26, American Mathematical Society, Providence, R.I., 1969.

\bibitem{Jenkins} J.A.\,Jenkins, \emph{Univalent functions and conformal
mapping}, Ergeb. Math. Grenzgeb. Neue Folge, vol.18, Reihe: Moderne Funktiontheorie,
Springer-Verlag, 1958.

\bibitem{Tamrazov1}  P.M.\,Tamrazov, \emph{A conformally metric theory of doubly-connected regions
and a generalized Blaschke product}, Dokl. Akad. Nauk SSSR 161(1965), pp.308--311(Russian).
English Translation: Soviet Math. Dokl. 6 (1965), pp.16432--16435.

\bibitem{Tamrazov2}   P.M.\,Tamrazov, \emph{On certain extremal problems in conformal mapping}, Mat. Sb. (N.S.),
73(115)(1967), pp.97--125(Russian).

\bibitem{DubininMZ} V.N.\,Dubinin, \emph{The Gr\"{o}tzsch and Teichm\"{u}ller extremal problems on a Riemann surface},
Mathematical Notes, 92(6)(2012), 103--110.

\bibitem{MMR} G.V.\,Milovanovi\'{c}, D.S.\,Mitrinovi\'{c}, Th.M.\,Rassias, \emph{Topics in polynomials: extremal problems,
inequalities, zeros}, World Scientific Publishing Co., Inc., Singapore, 1994.

\bibitem{BorweinErdelyi} P.\,Borwein, T.\,Erdelyi, \emph{Polynomials and polynomial inequalities}, Grad. Texts in
Math., 161, Springer-Verlag, New York, 1995.

\bibitem{RS} Q.I.\, Rahman, G.\,Schmeisser, \emph{Analytic theory of polynomials}, London Math.
Soc. Monogr. (N.S.), 26, The Clarendon Press, Oxford Univ.Press, 2002.

\bibitem{CourantHurwitz} A.\,Hurwitz, R.\,Courant, Vorlesungen \"{u}ber allgemeine Funktionentheorie und Elliptische Funktionen, Springer-Verlag, Berlin--New York, 1964.

\bibitem{DubKost} V.N.\,Dubinin and E.V.\,Kostyuchenko, \emph{The Teichm\"{u}ller extremal problem and distortion theorems in the theory of univalent functions}, Siberian Mathematical Journal,  40(2)(1999), pp.258--261.

\bibitem{DubSbornik} V.N.\,Dubinin, \emph{Inequalities for critical values of polynomials}, Sbornik:Mathematics, 197(8)(2006), pp.1167--1176.

\end{thebibliography}
\end{document}